\documentclass[11pt]{article}
\title{The First Ascent into the Incidence Algebra of the Fibonacci Cobweb Poset}
\author{Ewa Krot \\
\\Institute of Computer Science, Bia{\l}ystok University\\
PL-15-887 Bia{\l}ystok, ul.Sosnowa 64, POLAND\\
e-mail: ewakrot@wp.pl, ewakrot@ii.uwb.edu.pl}

\usepackage{amsmath,amsthm}

\chardef\bslash=`\\ % p. 424, TeXbook
\newtheorem{thm}{Theorem}[section]

\begin{document}
\maketitle
%%%%%%%%%%%%
\begin{abstract}
The explicite formulas for M\"{o}bius function and some other
important elements of the incidence algebra are delivered. For
that to do one uses Kwa\'sniewski's construction of his Fibonacci
cobweb poset in the plane grid coordinate system.

\end{abstract}
\section{Fibonacci cobweb poset}
The Fibonacci cobweb poset $P$ has been invented by
A.K.Kwa\'sniewski in \cite{4,3,10} for the purpose of finding
combinatorial interpretation of fibonomial coefficients and
eventually their reccurence relation. At first the partially
ordered set $P$ (Fibonacci cobweb poset) was there defined
 via Hasse diagram as follows: $P$ looks
like famous rabbits growth tree but it has the specific cobweb in
addition, i.e. it consists of levels labeled by Fibonacci numbers
(the $n$-th level consist of $F_{n}$ elements). Every element of
$n$-th level ($n \geq 0$) is in partial order relation with every
element of the $(n+1)$-th level but it's not with any element from
the level in which he lies ($n$-th level) except from it.

%\section{Plane grid coordinate system of $P$}
In \cite{4} A. K. Kwa\'sniewski defined cobweb poset $P$ as
infinite labeled digraph oriented upwards as follows: Let us label
vertices of $P$ by pairs of coordinates: $\langle i,j \rangle \in
{\bf N_{0}}\times {\bf N_{0}}$, where the second coordinate is the
number of level in which the element of $P$ lies (here it is the
$j$-th level) and the first one is the number of this element in
his level (from left to the right), here $i$. Following \cite{4}
we shall refer to $\Phi_{s}$ as to the set of vertices (elements)
of the $s$-th level, i.e.:
$$\Phi_{s}=\left\{\langle j,s \rangle ,\;\;1\leq j \leq F_{s}\right\},\;\;\;s\in{\bf N}\cup\{0\}.$$
For example $\Phi_{0}=\{\langle 1,0\rangle\},\;\Phi_{1}=\{\langle
1,1\rangle\},\; \Phi_{2}=\{\langle 1,2\rangle\},
\;\Phi_{3}=\{\langle 1,3\rangle, \langle
2,3\rangle\},\\\Phi_{4}=\{\langle 1,4\rangle, \langle 2,4\rangle,
\langle 3,4\rangle\},\;\Phi_{5}=\{\langle 1,5\rangle, \langle
2,5\rangle, \langle 3,5\rangle\,\langle 4,5\rangle,\langle
5,5\rangle\}$ ....

Then $P$ is a labeled graph $P=\left(V,E\right)$ where
$$V=\bigcup_{p\geq 0}\Phi_{p},\;\;\;E=\left\{\langle \,\langle j,p\rangle,\langle q,p+1
\rangle\,\rangle\right\},\;\;1\leq j\leq F_{p},\;\;1\leq q\leq
F_{p+1}.$$

 We can now define the partial order relation on $P$ as follows:
let\\ $x=\langle s,t\rangle, y=\langle u,v\rangle $ be elements of
cobweb poset $P$. Then
$$ ( x \leq y) \Longleftrightarrow
 [(t<v)\vee (t=v \wedge s=u)].$$

 \section{The Incidence Algebra ${\bf I(P)}$}
Let us recall that one defines the incidence algebra of a locally
finite partially ordered set $P$ as follows (see \cite{7,8}):
$$ {\bf I(P)}=\{f:{\bf P}\times {\bf P}\longrightarrow {\bf R};\;\;\;\; f(x,y)=0\;\;\;unless\;\;\; x\leq y\}.$$
The sum of two such functions $f$ and $g$ and multiplication by
scalar are defined as usual. The product $h=f\ast g$ is defined as
follows:
$$ h(x,y)=(f\ast g)(x,y)=\sum_{z\in {\bf P}:\;x\leq z\leq y} f(x,z)\cdot g(z,y).$$
It is immediately verified that this is an associative algebra
over the field of reals (or any associative ring) {\bf R}.

The incidence algebra has an identity element $\delta (x,y)$, the
Kronecker delta. Also the zeta function of ${\bf P}$ defined for
any poset by:
$$ \zeta (x,y)=\Big\{\begin{array}{l}1\;\;for\;\;\; x \leq y\\0\;\;otherwise\end{array}$$
 is an element of ${\bf I(P)}$. For the Fibonacci cobweb poset $\zeta$ was expressed by $\delta$ in
\cite{4} from where quote the result:
\begin{equation}\label{dzeta}
\zeta =\zeta_{1}-\zeta_{0}
\end{equation}
 where for $x,y \in {\bf N_{0}}$:
\begin{equation}
\zeta_{1}(x,y)=\sum_{k=0}^{\infty}\delta(x+k,y)
\end{equation}
\begin{equation}
\zeta_{0}(x,y)=\sum_{k \geq 0}\sum_{s \geq 0}\delta
(x,F_{s+1}+k)\sum_{r=1}^{F_{s}-k-1}\delta (k+F_{s+1}+r,y).
\end{equation}
 Using coordinates more convenient formula was given by the present author in
\cite{9}:
\begin{equation}
\zeta(x,y)=\zeta\left( \langle s,t\rangle ,\langle u,v \rangle
\right)=\delta (s,u)\delta (t,v)+\sum_{k=1}^{\infty}\delta(t+k,v).
\end{equation}
\section{The M\"{o}bius function and M\"{o}bius inversion formula
on $P$}
 The knowledge of $\zeta$ enables us to construct other
typical elements of incidence algebra  of ${\bf P}$. The one of
them is M\"{o}bius function indispensable in numerous inversion
type formulas of countless applications \cite{7,8}. Of course the
$\zeta$ function of a locally finite partially ordered set is
invertible in incidence algebra and its inversion is the famous
M\"{o}bius function $\mu$ i.e.:
$$\zeta \ast \mu=\mu \ast \zeta=\delta.$$The M\"{o}bius function
$\mu$ of Fibonacci cobweb poset ${\bf P}$ was presented  for the
first time by the present author in \cite{2}. It was recovered
just by the use of the recurrence formula for M\"{o}bius function
of locally finite partially ordered set  ${\bf I(P)}$ (see
\cite{7}):
\begin{equation}
 \left\{\begin{array}{l}\mu(x,x)=1\;\;\;\;for\;\; all \;\;x\in{\bf
 P}\quad \quad \quad \\\ \\
\mu (x,y)=-\sum_{x\leq z<y}\mu (x,z)\end{array}\right.
\end{equation}

The grid coordinates system definition of $P$ \cite{4} allowed the
present author to derive an explicit formula for M\"{o}bius
function of cobweb poset $P$, \cite{9}.
 Namely for $x=\langle s,t\rangle,\;\;y=\langle u,v\rangle,\,\,1\leq s\leq
F_{t},\;\;1\leq u\leq F_{v},\;t,v\in {\bf N_{0}}$ we have
\begin{multline}\label{mobius2}
\mu(x,y)=\mu\left( \langle s,t\rangle ,\langle u,v\rangle
\right)=\\
\qquad \quad \;
=\delta(t,v)\delta(s,u)-\delta(t+1,v)+\sum_{k=2}^{\infty}\delta(t+k,v)(-1)^{k}
\prod_{i=t+1}^{v-1}(F_{i}-1)
\end{multline}
where $\delta$ is the usual Kronecker delta
$$ \delta (x,y)=\Big\{\begin{array}{l}1\;\;\;\;\; x=y\\0\;\;\;\;\;x\neq y\end{array}.$$

The formula (\ref{mobius2}) enables us to formulate the following
theorem (see \cite{7}):
\begin{thm}{\bf (M\"{o}bius Inversion Formula for $P$)}\\
Let $f(x)=f(\langle s,t\rangle)$ be a real valued function,
defined for $x=\langle s,t\rangle$ ranging in cobweb poset $P$.
Let an element $p=\langle p_{1},p_{2}\rangle$ exist with the
property that $f(x)=0$ unless $x\geq p$.

Suppose that $$g(x)=\sum_{\left\{ y \in P :\; y\leq x
\right\}}f(y).$$ Then
$$f(x)=\sum_{\{y\in P:\;y\leq x\}}g(y)\mu(y,x).$$
Hence using coordinates of $x,y$ in $P$ i.e. $x=\langle
s,t\rangle,\;y=\langle u,v\rangle$ if
$$g(\langle s,t\rangle)=\sum_{v=0}^{t-1}\sum_{u=1}^{F_{v}}\left(f(\langle u,v\rangle )\right)+f(\langle s,t\rangle)$$
then we have
\begin{multline}
f(\langle s,t\rangle)=\sum_{v\geq 0}\sum_{u=1}^{F_{v}}g(\langle
u,v\rangle)\mu(\langle s,t\rangle ,\langle
u,v\rangle)=\\
=\sum_{v\geq 0}\sum_{u=1}^{F_{v}}g(\langle
u,v\rangle)\left[\delta(v,t)\delta(u,s)-\delta(v+1,t)+\sum_{k=2}^{\infty}\delta(v+k,t)(-1)^{k}
\prod_{i=v+1}^{t-1}(F_{i}-1)\right].
\end{multline}
\end{thm}
\section{Examples of Other Elements of ${\bf I(P)}$}
 The knowledge of $\zeta$ enables one to construct the typical
elements of incidence algebra perfectly suitable for calculating
number of chains, of maximal chains etc. in {\bf finite}
sub-posets of ${\bf P}$. We shall then now deliver some of them in
explicit form.
% \\{\bf CARDINALITY OF SEGMENT
\renewcommand{\labelenumi}{(\arabic{enumi})}
\begin{enumerate}
\item The function $\zeta ^{2}=\zeta \ast \zeta$ counts the number
of elements in the segment $\left[ x,y\right]$
\\(where $x=\langle s,t \rangle ,
y=\langle u,v \rangle$), i.e.:
\begin{multline*}
\zeta^{2}(x,y)=\left(\zeta \ast \zeta\right) (x,y)=\\\quad
\quad=\sum_{x \leq z \leq y}\zeta (x,z)\cdot\zeta(z, y)=\\=\sum_{x
\leq z \leq y}1=\mathrm{card} \left[ x, y\right]\quad \quad \quad
\;\;\,
\end{multline*}
Therefore for $x,y\in {\bf P}$ as above, we have:
$$card \left[
x, y\right]=\left(\sum_{i=t+1}^{v-1}\sum_{j=1}^{F_{i}}1\right)+2$$
% {\bf NUMBER OF CHAINS IN A SEGMENT}
\item  For any incidence algebra the function $\eta$ is defined as
follows: $$\eta(x,y)=(\zeta-\delta)(x,y)=\left\{
\begin{array}{lll}
1&&x<y\\
0&&otherwise \end{array}\right.$$ Hence $\eta ^{k}(x,y),\;\; (k\in
{\bf N})$ counts the number of chains of lenght $k$, (with $(k+1)$
elements) from $x$ to $y$.

 The
corresponding function for $x,y$ being elements of Fibonacci
cobweb poset ${\bf P}$, ($x=\langle s,t \rangle , y=\langle u,v
\rangle$) is then given by formula:
$$ \eta(\langle s,t\rangle , \langle u,v\rangle)=\sum_{k=1}^{\infty}\delta(t+k,v).$$
\item  Now let
$$\mathcal{C}(x,y)=(2\delta-\zeta)(x,y)=\left\{\begin{array}{lll}
1&&x=y\\
-1&&x<y\\
0&&otherwise \end{array} \right.$$ For elements of ${\bf P}$ we
have:
$$\mathcal{C}(\langle s,t\rangle , \langle u,v\rangle)=
\delta(t,v)\delta(s,u)-\sum_{k=1}^{\infty}\delta(t+k,v)$$

The inverse function $\mathcal{C}^{-1}(x,y)$ counts the number of
all chains from $x$ to $y$. From the recurrence formula one infers
that
$$\left\{\begin{array}{l}
\mathcal{C}^{-1}(x,x)=\frac{1}{\mathcal{C}(x,x)}\\
\\
\mathcal{C}^{-1}(x,y)=-\frac{1}{\mathcal{C}(x,x)}\sum_{x< z\leq
y}\mathcal{C}(x,z)\cdot\mathcal{C}^{-1}(z,y) \end{array}\right.$$
%{\small
%\begin{center}
%{\bf NUMBER OF MAXIMAL CHAINS IN A SEGMENT}
%\end{center}
\item For any incidence algebra the function $\chi$ is defined as
follows:
$$\chi (x,y)=\left\{\begin{array}{lll}
1&&y\;\; covers\;\;x\\
0&&otherwise. \end{array}\right. $$ Therefore $\chi ^{k}(x,y),\;\;
(k\in {\bf N})$ counts the number of maximal chains of lenght $k$,
(with $(k+1)$ elements) from $x$ to $y$.

 The
corresponding function for $x,y$ being elements of Fibonacci
cobweb poset ${\bf P}$ ($x=\langle s,t \rangle , y=\langle u,v
\rangle$) is expressed by formula:
$$ \chi(\langle s,t\rangle , \langle u,v\rangle)=\delta(t+1,v).$$

Finally let

$$\mathcal{M}(x,y)=(\delta -\chi)(x,y)=\left\{
\begin{array}{ll}
1&x=y\\
-1&y\;\;covers\;\;x\\
0&otherwise \end{array}\right.$$ For elements of ${\bf P}$ we
have:
$$\mathcal{M}(\langle s,t\rangle ,\langle u,v\rangle)=\delta(t,v)\delta(s,u)-\delta(t+1,v).$$
Then the  inverse function of $\mathcal{M}$:
$$\mathcal{M}^{-1}=\frac{\delta}{\delta -\chi}=\delta +\chi +\chi
^{2}+\chi^{3}+\ldots$$ counts the number of all maximal chains
from $x$ to $y$.
\end{enumerate}
{\bf Closing Remark}\\
As the poset $P$ is not of binomial type \cite{11} the further
study of the $I(P)$ algebra looks promising as for as obstacles
are concerned.
\\
\\
{\bf Acknowledgements}\\ The author is indebted to professor
A.K.Kwa\'sniewski also for suggesting the title if this note.

 AMS Classification numbers: 11C08, 11B37, 47B47
 \end{document}